\documentclass[12pt]{article}

\usepackage{amsmath}
\allowdisplaybreaks[4]

\usepackage[all]{xy}
\usepackage{amssymb}
\usepackage{amsthm}
\usepackage{hyperref}
\hypersetup{colorlinks=true,linkcolor=blue,citecolor=red}
\usepackage{amsmath}
\usepackage{amscd}
\usepackage{verbatim}
\usepackage{eurosym}
\usepackage{float}
\usepackage{color}
\usepackage{dcolumn}
\usepackage[mathscr]{eucal}
\usepackage[all]{xy}
\usepackage{hyperref}
\usepackage{mathrsfs}
\usepackage{amsmath}
\usepackage{amssymb}
\usepackage{amsfonts,ifpdf}
\usepackage{graphicx}
\usepackage{times}
\usepackage{float}
\usepackage{epstopdf}
\usepackage{cite}

\newtheorem{theo}{Theorem}[section]

\newtheorem{lem}[theo]{Lemma}

\newtheorem{con}[theo]{Conjecture}

\newtheorem{pf}{proof}[section]
\theoremstyle{remark}

\makeatletter \@addtoreset{equation}{section} \makeatother
\makeindex 
\def\qed{\hfill \rule{4pt}{7pt}}

\begin{document}

\begin{center}
{\Large\bf Semi-invariants of Binary Forms Pertaining to\\[5pt]
a Unimodality Theorem of Reiner and Stanton
}\\ [7pt]
\end{center}

\vskip 3mm

\begin{center}
William Y.C. Chen$^{1}$ and Ivy D.D. Jia$^{2}$\\[8pt]
$^{1,2}$Chern Institute of Mathematics
\\
Nankai University\\
Tianjin 300071, P. R. China\\[3pt]
and\\[3pt]
$^{1}$Center for Applied Mathematics\\
Tianjin University\\
Tianjin 300072, P. R. China\\[15pt]

Emails: $^{1}$chenyc@tju.edu.cn, $^{2}$jiadandan@mail.nankai.edu.cn
\\[15pt]
\textit{Dedicated to the Memory of Professor S. S. Chern } \\[15pt]

\end{center}

\vskip 3mm

\begin{abstract}

The symmetric difference of the $q$-binomial coefficients $F_{n,k}(q)={n+k\brack k}-q^{n}{n+k-2\brack k-2}$ was introduced by Reiner and Stanton.
They proved that $F_{n,k}(q)$ is symmetric and unimodal for $k \geq 2$ and $n$ even by using the representation theory for Lie algebras.
Based on Sylvester's proof of the unimodality of the Gaussian coefficients, as conjectured by Cayley, we find
an interpretation of the unimodality of $F_{n,k}(q)$ in terms of
semi-invariants.
In the spirit of the strict unimodality of the Gaussian coefficients due to Pak and Panova, we prove the strict unimodality of the symmetric
difference
$G_{n,k,r}(q)={n+k\brack k}-q^{nr/2}{n+k-r\brack k-r}$, except for the two terms at both ends, where $n,r\geq8$, $k\geq r$ and at least one of $n$ and $r$ is even.

\vskip 6pt

\\ [7pt]
\end{abstract}

\section{Introduction}

The story begins with  the Gaussian coefficients called by Rota with no particular reasons (private conversation, see also \cite{GR69}), or sometimes the Gaussian polynomials, or often under the name of the $q$-binomial coefficients, as given by
\begin{equation*}
{n\brack k}=\frac{(1-q^{n})(1-q^{n-1})\cdots(1-q^{n-k+1})}{(1-q^{k})(1-q^{k-1})\cdots(1-q)},
\end{equation*}
where $0\leq k \leq n$.
The Gaussian coefficients are polynomials in $q$ and they enjoy the fundamental property that their coefficients are symmetric and unimodal, as conjectured by Cayley \cite{Cayley56} in 1856 and  confirmed by Sylvester \cite{Sylvester78} in 1878, who had even believed that settling the conjecture of Cayley
was a task that lay outside the human power, see also, Pak and Panova \cite{PP13a,PP13b}.
Ever since a great deal of work has been done in this vein, see, for example,
\cite{Bressoud92,Dhand14,DK17,DK18,GOS92,Kirillov92, O'Hara90,PP14,PPV14,PP17,PR86,Proctor82,Stanley80,Stanley89,Stanton90, Zanello15, Zeilberger89}, to mention only a few.
In particular, O'Hara \cite{O'Hara90} found a combinatorial proof,
 Zeilberger \cite{Zeilberger89}
came up with an identity, known as the KOH theorem, which serves the purpose
of justifying the unimodality.

The unimodality is not only associated with the Gaussian coefficients, it can also be said about certain differences of the Gaussian coefficients.
Employing the representation theory for Lie algebras, Reiner and Stanton \cite{RS98}  established the unimodality of
\begin{align}
{n\brack k}-{n\brack k-1},
\end{align}
where $k\geq 1$, $n$ is odd and $2k\leq n+1$.
They also showed that
\begin{equation}\label{Fnk}
F_{n,k}(q)={n+k\brack k}-q^{n}{n+k-2\brack k-2}
\end{equation}
is symmetric and unimodal when $k\geq 2$ and $n$ is even.
Furthermore, Reiner and Stanton conjectured the unimodality of
\begin{align}\label{k-2}
{n-1\brack k}-q^{n-2rk+1+4(r-1)}{n-1+4(r-1)\brack k-2},
\end{align}
where $n$ is odd, $k\geq 2$, $r\geq 1$ and  $n\geq 2rk-4r+3$.
This conjecture is still open.

The above symmetric differences \eqref{k-2} are called the strange symmetric differences by
Stanley and Zanello \cite{SZ20}.
They extended the above conjecture  to a broader framework, namely, for each $k\geq 5$, the polynomials
\begin{equation} \label{fkmb}
f_{k,m,b}(q)={m\brack k}-q^{\frac{k(m-b)}{2}+b-2k+2}{b\brack k-2}
\end{equation}
are nonnegative and unimodal subject to certain conditions.
They verified their conjecture for $2\leq k\leq 5$ by means of the KOH theorem of Zeilberger \cite{Zeilberger89}.
Notice that when $b=m-2$, \eqref{fkmb} reduces to \eqref{Fnk}.

It is worth mentioning that Bergeron \cite{Bergeron16} investigated the symmetric differences  of the following form:
\begin{align}\label{abcd}
{b+c\brack b}-{a+d\brack d},
\end{align}
where $a,b,c,d$ are positive integers, with $a$ being the smallest and $ad=bc$.
It was conjectured that such symmetric differences are polynomials in $q$ with nonnegative coefficients.
Zanello \cite{Zanello18} further conjectured the unimodality of \eqref{abcd}, and proved the unimodality for $a\leq3$ and $b,c\geq4$, resorting to the KOH theorem.

A notable progress on the unimodality of the Gaussian coefficients was achieved by Pak and Panova \cite{PP13a,PP13b}.
They showed that the Gaussian coefficients are strictly unimodal except for a few cases.
Recall that a polynomial $f(x)=a_0 + a_1x + a_2x^2 + \cdots + a_n x^n$ with positive coefficients is said to be unimodal if there is an index $m$ such that
\begin{equation} \label{u-1}
 a_0\leq a_1\leq \cdots \leq a_m \geq a_{m+1} \geq  \cdots \geq a_n .
\end{equation}
We say that $f(x)$ is strictly unimodal if \eqref{u-1} is to be replaced by
\begin{equation} \label{u-2}
 a_0<  a_1<  \cdots <  a_m > a_{m+1} >  \cdots > a_n .
\end{equation}

The humble goal of this paper is to present an interpretation of the unimodality of $F_{n,k}(q)$ as  in \eqref{Fnk} in terms of semi-invariants of binary forms, following the original aspiration of Sylvester.
Once the connection to semi-invariants is at our disposal, we carry on to prove the strict unimodality of the symmetric difference
\begin{equation}\label{Gnkr}
G_{n,k,r}(q)={n+k\brack k}-q^{\frac{nr}{2}}{n+k-r\brack k-r},
\end{equation}
except for the two terms at both ends, when $n,r\geq8$, $k\geq r$ and at least one of $n$ and $r$ is even.

It is our hope that after a sound rest for more than a century, it might be the time
for the once shining invariant theory of binary forms to shed light on the ongoing study of the
theory of partitions, presumably
under a banner bearing the name of enumerative invariant theory.

\section{Semi-invariants}

A binary form of degree $n$, or a binary $n$-form, is a homogeneous polynomial in $x$ and $y$,
\begin{equation} \label{F-1}
f(x,y)= a_0x^n+ \binom n1 a_{1} x^{n-1}y+{n\choose 2} a_{2} x^{n-2}y^2+\cdots+a_ny^n.
\end{equation}
A semi-invariant or a source of a covariant, called by Cayley   \cite{Cayley81}, or a differentiant
  called by Sylvester \cite{Sylvester78}, of the binary form $f(x,y)$, is
a polynomial $I(a_0, a_1, \ldots, a_n)$ in $a_0, a_1, \ldots, a_n$
with rational coefficients
that is invariant under the transformation:
\begin{equation} \label{s-trans}
\left( \begin{matrix}
        x \\[3pt]
        y
        \end{matrix}
        \right) =  \left(
     \begin{matrix}
1& h\\[3pt]
0& 1
\end{matrix}
\right)
\left( \begin{matrix}
        x' \\[3pt]
        y'
        \end{matrix}
        \right)
,
\end{equation}
that is, $x= x' + h y'$ and $y=y'$.
The matrix in \eqref{s-trans} is called a shear matrix.
Notice that Hilbert \cite{Hilbert93} used the Greek letter $\mu$ in place
of $h$ in the  shear matrix in \eqref{s-trans}.
Here we prefer the above notation of Sylvester \cite{Sylvester78}
for the reason that
the  Greek letters $\lambda, \mu$ are now
often reserved for integer partitions. Indeed, the transformation \eqref{s-trans} for binary forms goes back to
 Lagrange, see \cite{Grosshans03,Salmon85}. More precisely,
suppose that the binary form \eqref{F-1} becomes
\begin{equation} \label{F-2}
f'(x',y')= a'_0x'^n+ \binom n1 a'_{1} x'^{n-1}y'+{n\choose 2} a'_{2} x'^{n-2}y'^2+\cdots+a'_ny'^n
\end{equation}
under the transformation \eqref{s-trans}, where, for $0\leq i\leq n$,
\begin{equation}
a_i'=a_i+{i\choose1} a_{i-1}h+{i\choose2} a_{i-2}h^2+\cdots+a_0 h^i.
\end{equation}
Then we say that
a polynomial $I(a_0, a_1, \ldots, a_n)$ with rational coefficients
is a semi-invariant
of the binary form \eqref{F-1} if
\begin{equation}\label{II}
I(a_0, a_1, \ldots, a_n) = I(a'_0, a'_1, \ldots, a'_n)
\end{equation}
for any shear transformation in \eqref{s-trans}.

For example, for the quadratic form
\begin{equation}
f(x,y)=a_0x^2+2a_1xy+a_2y^2,
\end{equation}
it can be easily checked that the polynomial
\begin{equation}\label{I-a}
 {I}(a_0,a_1,a_2)=a_0a_2-a_1^2
\end{equation}
is a semi-invariant. In fact, it is the discriminant, which is more than a semi-invariant in
the sense that it is an invariant with respect to a more
general transformation, see, for example, \cite{Grosshans03,Sturmfels08}.
Below is a semi-invariant
but not an invariant:
\begin{equation}\label{J-a}
J(a_0,a_1,a_2)= a_0^2a_2-a_0a_1^2.
\end{equation}

As far as this paper is concerned, a semi-invariant can be viewed as
nothing but a polynomial $I(a_0,a_1, \ldots, a_n)$
 that satisfies a
 partial differential equation. Needless to say, the coefficients
  $a_0, a_1, \ldots, a_n$
 are perceived as variables. First, the degree of a monomial
\begin{equation}
a^\nu=a_0^{\nu_0}a_1^{\nu_1}\cdots a_n^{\nu_n}
\end{equation}
is meant to be the total degree $\nu_0 + \nu_1+\cdots +\nu_n$.
Besides,
the weight of the  monomial $a^\nu$ is defined by
\begin{align}
\nu_1+2 \nu_2+\cdots+n \nu_n.
\end{align}
For the purpose of this paper,
we shall only be concerned with homogeneous semi-invariants with respect to the
 degree and the weight.
 For example, the semi-invariant
$I(a_0, a_1, a_2)$ in \eqref{I-a} has degree two and weight two,
whereas the semi-invariant $J(a_0, a_1, a_2)$ in \eqref{J-a}
has degree three and weight two.

We are now led to define  $Q_n(k,m)$ as the
vector space of polynomials in $a_0, a_1, \ldots, a_n$
over the rational numbers that are homogeneous of degree $k$ and weight $m$.
Taking the dimension of $Q_n(k,m)$ into account,
the notion of partitions immediately comes into play.
The monomials
 in $Q_n(k,m)$ are in one-to-one
 correspondence with partitions of $m$ into $k$ parts with zero parts allowed and each part
 not exceeding $n$.
   Recall that a partition of $m$ with $k$ parts such that
 no parts exceed  $n$
 can be represented as
 \begin{equation}
 \lambda=0^{\nu_0} 1^{\nu_1} \cdots n^{\nu_n},
 \end{equation}
 where $\nu_i$ signifies the number of occurrences of
 the part $i$ in $\lambda$ and $\nu_0+\nu_1+\cdots + \nu_n=k$.
 We say that such a partition is contained in a $k\times n$ rectangle.
 Clearly, $\lambda$ is a partition of $\nu_1 +2\nu_2+
 \cdots + n\nu_n$, which is the weight of
 the monomial $a^\nu$. This basic fact
 brings us right to the stage of the
 Gaussian coefficients which admit the following partition
 interpretation.

Let $p(k,n,m)$ denote the number of partitions of $m$ contained in a $k\times n$ rectangle, then we have
\begin{equation} \label{nkk}
{n+k\brack k}=\sum_{m=0}^{nk}p(k,n,m)q^m,
\end{equation}
see \cite{Andrews98,Stanley12}. The above
polynomial is symmetric, that is,
for $0\leq m\leq nk/2$,
\begin{equation}\label{symmetry}
p(k,n,m)=p(k,n,nk-m).
\end{equation}
For $0\leq m\leq nk/2$,  Sylvester found an interpretation of the difference $p(k,n,m)-p(k,n,m-1)$
in terms of semi-invariants, which demystifies the unimodality
of the Gaussian coefficients.

The following characterization is crucial for the route from the vector space $Q_n(k,m)$ to semi-invariants, see Cayley \cite{Cayley56} or
Hilbert \cite{Hilbert93}.

\begin{theo} \label{DI}
Let
\begin{equation} \label{D}
 {D} = a_0\frac{\partial}{\partial a_1}+2a_1\frac{\partial}{\partial a_2}
+3a_2\frac{\partial}{\partial a_3}+\cdots+na_{n-1}\frac{\partial}{\partial a_n}.
\end{equation}
 A polynomial $I(a_0,a_1,\ldots,a_n)$ in $Q_n(k,m)$
 is a semi-invariant of the binary  form \eqref{F-1} if and only if $ {D}I(a_0,a_1,\ldots,a_n)=0$.
\end{theo}

 Taking the semi-invariant $J(a_0,a_1, a_2)$ in \eqref{J-a} as
 an example, it can be seen that
\begin{equation}
{D} J(a_0,a_1, a_2)=\left(a_0\frac{\partial}{\partial a_1}+2a_1\frac{\partial}{\partial a_2}\right)(a_0^2a_2-a_0a_1^2)
\end{equation}
is identically zero.

Let $S_n(k,m)$
denote the set of semi-invariants of degree $k$ and
weight $m$, that is,
\begin{equation} \label{Vkp}
S_n(k,m) = \{I \in Q_n(k,m)  \mid D (I) = 0 \},
\end{equation}
where $D$ is the linear operator defined in \eqref{D}.
Bear in mind that $S_n(k,m)$ forms a
vector space over $\mathbb{Q}$.
By the number of
semi-invariants of degree $k$ and weight $m$ of
a binary $n$-form, we really
mean the dimension of the vector space $S_n(k,m)$.
For example, $\dim S_4(4,6)=2$. Below are two
linearly independent semi-invariants of degree four
and weight six of a binary $4$-form:
\begin{align}\label{I1I2}
I_1&=3a_1^2 a_2^2-4a_1^3 a_3-2a_0 a_1 a_2 a_3
+3a_0^2 a_3^2+4a_0a_1^2a_4-4a_0^2a_2a_4,\\[6pt]
I_2&=a_0 a_2^3-2a_0 a_1 a_2 a_3+a_0^2 a_3^2
+a_0 a_1^2 a_4-a_0^2 a_2 a_4.\label{I1I2-b}
\end{align}
We now come to the remarkable discovery of
Sylvester \cite{Sylvester78}.

\begin{theo} \label{weishu}
For $n,k\geq 0$ and $0\leq m\leq nk/2$, the number of semi-invariants of a binary $n$-form of degree $k$ and weight $m$ equals
\begin{equation}
\delta(k,n,m) = p(k,n,m)-p(k,n,m-1),
\end{equation}
with the convention that $p(k,n,-1)=0$.
\end{theo}

For example, for the quadratic form
\begin{equation}
f(x,y)=a_0x^2+2a_1xy+a_2y^2,
\end{equation}
there is exactly
one semi-invariant of degree three and weight two, that is, \begin{equation}
J(a_0, a_1, a_2)= a_0^2a_2-a_0a_1^2.
\end{equation}
On the other hand, we have
 $p(3,2,2)=2$ and  $p(3,2,1)=1$, in accordance with Theorem \ref{weishu}.

Up to now, we are sufficiently equipped to move on to the next section to
explore symmetric differences of the Gaussian coefficients by means of semi-invariants.

\section{Symmetric Differences of  the Gaussian Coefficients}

The first objective of this section is to give a
semi-invariant interpretation of the unimodality theorem of Reiner and Stanton on $F_{n,k}(q)$ as in \eqref{Fnk}. Then we proceed to
prove the strict unimodality of $G_{n,k,r}(q)$ as in \eqref{Gnkr}.

\begin{lem} \label{l-even}
If $n$ is even, then there is exactly one semi-invariant of a binary $n$-form of degree two and weight $n$.
\end{lem}

\noindent
{\it Proof.}
It is not hard to see that for $m\leq n$,
\begin{equation} \label{p2}
p(2,n,m)=\left\lfloor { m+2 \over 2} \right\rfloor,
\end{equation}
see Stanley and Zanello  \cite{SZ20, Zanello18}.
Since $n$ is even, by Theorem \ref{weishu}  we see that
there is only one semi-invariant of degree two and weight $n$.
\qed

For example, there is only one semi-invariant of a binary $4$-form of degree two and weight four, that is,
\begin{equation}
J=3a_2^2-4a_1 a_3 +a_0 a_4.
\end{equation}

The argument in the proof of Lemma \ref{l-even}
 indicates that when $n$ is odd,
there are no semi-invariants of a binary $n$-form of degree two and weight $n$.
The following relation is essentially a consequence of the fact that the set
of semi-invariants forms a ring.
To be precise, if $I$ and $J$ are two semi-invariants of a binary $n$-form, then so are $I+J$ and $IJ$.
In fact, the ring property of semi-invariants
is a consequence of the operator characterization as stated in
Theorem \ref{DI}.

\begin{theo}\label{uni}
If $k\geq2$, $m\geq n$ and $n$ is even, then the number of semi-invariants of a binary $n$-form of degree $k$ and weight $m$ is at least the number of semi-invariants of degree $k-2$ and weight $m-n$, that is,
\begin{equation} \label{dk2}
\delta(k,n,m) \geq \delta(k-2, n, m-n).
\end{equation}
\end{theo}

\noindent
{\it Proof.} Assume that there are $t$ linearly independent
semi-invariants of degree $k-2$ and weight $m-n$,
say, $I_1, I_2, \ldots, I_t$, where the variables $a_0, a_1, \ldots, a_n$
are suppressed. We wish to show that there are at least
$t$ semi-invariants of degree $k$ and weight $m$. 
By Lemma \ref{l-even},
 we may assume that  $J$ is a semi-invariant
 of a binary $n$-form of degree two  and weight $n$.
Thanks to the ring structure of semi-invariants,
 we see that $J I_1, JI_2, \ldots , JI_t$ are semi-invariants
 of degree $k$ and weight $m$.
To complete the proof, one only needs to realize that
$JI_1,\,JI_2,\,\ldots,JI_{t}$  are linearly independent, which is by any means
 a plain fact.
\qed

Next we demonstrate that Reiner and Stanton's unimodality theorem for $F_{n,k}(q)$ is   immediate  from Theorem \ref{uni}.
Recall that
\begin{equation}\label{Fnk-b}
F_{n,k}(q)={n+k\brack k}-q^{n}{n+k-2\brack k-2},
\end{equation}
where $k\geq 2$ and $n$ is even.

While keeping the symmetry of $F_{n,k}(q)$ in mind,
to confirm the unimodality, let
\begin{equation}
F_{n,k}(q)=\sum_{m=0}^{nk}f_mq^m,
\end{equation}
so that
\begin{equation}
f_m=p(k,n,m)-p(k-2,n,m-n),
\end{equation}
with the convention that $p(k-2,n, j)=0$ if $j$ is negative.
For $0\leq  m\leq nk/2$, in light of Theorem \ref{weishu},
the inequality \eqref{dk2} yields
\begin{equation}
p(k,n,m)-p(k,n,m-1)\geq p(k-2,n,m-n)-p(k-2,n,m-n-1),
\end{equation}
which can be recast as
\begin{equation}
p(k,n,m)-p(k-2,n,m-n)\geq p(k,n,m-1)-p(k-2,n,m-n-1).
\end{equation}
But this is exactly $f_m\geq f_{m-1}$.

Notice that $F_{n,k}(q)$ is not always strictly unimodal, except for the two terms at both ends.
For example, for $F_{5,9}(q)$ and $F_{14,5}(q)$, the maximal
coefficients occur  more than twice in the middle. The following
conjecture is supported by numerical evidence.

 \begin{con} For $n\geq 8$ and $k\geq 15$,
 $F_{n,k}(q)$ is strictly unimodal, except for the two terms at both ends.
 \end{con}

One might expect to push forward along this direction to
 tackle the strict unimodality of $F_{n,k}(q)$.
But it is not clear as to  how this can be pursued.
Nevertheless, we are given a chance to establish
     the strict unimodality of $G_{n,k,r}(q)$, with the two terms at the beginning and at the end being excluded,
     to be precise.
    To this end, we need two special
  semi-invariants, as ensured by the following lemma.

\begin{lem}\label{nr8}
If $n,r\geq8$ and at least one of $n$ and $r$ is even, then there are at least two linearly independent
 semi-invariants of a binary $n$-form of degree $r$ and weight $nr/2$, that is,
 \begin{equation} \label{d-r}
\delta(r,n,nr/2) \geq 2.
 \end{equation}
\end{lem}

In a more general setting,
 Pak and Panova \cite{PP13a,PP13b} obtained the strict unimodality of the Gaussian coefficients.

\begin{theo} \label{1}
For all $n,k\geq 8$ and $2\leq m\leq nk/2$, we have $\delta(k,n,m)\geq 1$.
\end{theo}

It should be pointed out that the following proof of Lemma \ref{nr8} is reminiscent of the reduction
argument of Pak and Panova.

\noindent
{\it Proof of Lemma \ref{nr8}.}
For $8 \leq n , r < 16$,  it is easily verified that \eqref{d-r} holds.
We now turn to the case
$8 \leq n < 16$ and $r\geq 16$.
Write $r=8s+t$, where $s\geq 1$ and $8\leq t< 16$. As is known for the first case,
we have $\delta(8,n,4n)\geq 2$.
This means that there exists a semi-invariant
$I$, with the variables being suppressed, of a binary $n$-form of degree eight and weight $4n$.
Again, the argument for the first case implies
$\delta(t,n,nt/2)\geq 2$. This enables us
to find two linearly independent semi-invariants $J_1, J_2$ of a binary $n$-form of degree $t$ and weight $nt/2$. Observe that
$I^sJ_1, I^sJ_2$ are linearly
independent
semi-invariants of a binary $n$-form of degree $r$ and weight $nr/2$,
from which we obtain that $\delta(r,n,nr/2)\geq 2$ for $8\leq n<16$ and $r\geq 16$.
Together with the first case $8\leq n,r<16$, we see that $\delta(r,n,nr/2)\geq2$
for $ 8\leq n < 16$ and $r\geq 8$.
   Since
\begin{equation}
p(r,n,nr/2)=p(n,r,nr/2),
\end{equation}
$\delta(r,n,nr/2)\geq2$ also holds
for $n\geq 8$ and  $8\leq r < 16$.

 Now we are left only with the case $n\geq 8$ and $r\geq 16$. For $n\geq 8$, we  write $r=8s+t$,  where $s\geq 1$ and $8\leq t< 16$. Mimicking the above
 reasoning for  $8 \leq n < 16$ and $r\geq 16$, we may reduce this case
 ($n\geq 8$ and $r\geq 16$) back to the case $8\leq n<16$ and $r\geq 16$.
 In summary, we conclude that $\delta(r,n,nr/2)\geq 2$ holds for all $n,r\geq 8$.
\qed

\begin{theo}\label{strict}
If $n,r\geq8$, $k\geq r$, $m\geq nr/2$, and at least one of $n$ and $r$ is even, then the number of semi-invariants of a binary $n$-form of degree $k$ and weight $m$ is greater than the number of semi-invariants of degree $k-r$ and weight $m-nr/2$, that is,
\begin{equation}\label{nr/2}
\delta(k,n,m)> \delta(k-r,n,m-nr/2).
\end{equation}
\end{theo}

Let us see how the above relation yields the strict unimodality of $G_{n,k,r}(q)$.

\begin{theo} For $n,r\geq 8$, $k\geq r$ and at least one of $n$ and $r$ is even,
\begin{equation}\label{Gnkr-b}
G_{n,k,r}(q)={n+k\brack k}-q^{\frac{nr}{2}}{n+k-r\brack k-r}
\end{equation}
is symmetric and strictly unimodal, except for the
two terms at both ends.
\end{theo}

\noindent
{\it Proof.}
The symmetry of $G_{n,k,r}(q)$ can be easily verified.
Let
\begin{equation}
G_{n,k,r}(q)=\sum_{m=0}^{nk}g_mq^m,
\end{equation}
and so
\begin{equation}
g_m=p(k,n,m)-p(k-r,n,m-nr/2),
\end{equation}
where we assume that $p(k-r,n, j)=0$ whenever $j$ is negative.
Notice that $g_0=g_1=1$ and $g_{nk-1}=g_{nk}=1$.
For $2 \leq m < nr/2$, we have $g_m=p(k,n,m)$ since $p(k-r,n,m-nr/2)=0$. Thus, by Theorem
\ref{1} we see that $g_m-g_{m-1}=\delta(k,n,m)$, which is greater than or equal to one. 
For $nr/2 \leq m \leq nk/2$, in virtue of Theorem \ref{weishu},
the inequality \eqref{nr/2} takes the form
\begin{equation}
p(k,n,m)-p(k,n,m-1)>  p(k-r,n,m-nr/2)-p(k-r,n,m-nr/2-1),
\end{equation}
which can be reformulated as $g_{m} >  g_{m-1}$. Thus we have shown that
$G_{n,k,r}(q)$ is symmetric and strictly unimodal. \qed

For example, we have
\begin{align*}
G_{8,14,10}(q)=&1+ q + 2 q^2 + 3 q^3 + 5 q^4 + 7 q^5 + 11 q^6 + 15 q^7 \\[5pt]
&\quad+ \cdots + 8310 q^{53} + 8408 q^{54} + 8450 q^{55} + 8479 q^{56} \\[5pt]
&\quad +
 8450 q^{57} + 8408 q^{58} + 8310 q^{59}+ \cdots
 + 11 q^{106} \\[5pt]
&\quad + 7 q^{107} + 5 q^{108} + 3 q^{109} +
 2 q^{110} + q^{111} + q^{112}.
\end{align*}

To present the proof of Theorem \ref{strict},
it is necessary to define the leading term of a polynomial $f(a_0, a_1, \ldots, a_n)$.
First, we write a monomial in the form $a^\nu=a_0^{\nu_0} a_{1} ^{\nu_{1}} \cdots a_n^{\nu_n}$.
Then we order the monomials of degree $k$ and weight $m$ according to the anti-lexicographic order
of their exponents. Clearly, this order extends to the set of all monomials
 in $a_0, a_1, \ldots, a_n$. For example, for $k=4$, $n=4$ and $m=6$, we have
\begin{align*}
a_1^2 a_2^2>a_0 a_2^3>a_1^3 a_3>a_0 a_1 a_2 a_3>a_0^2 a_3^2>a_0 a_1^2 a_4>a_0^2 a_2 a_4.
\end{align*}
The leading term of a semi-invariant $I(a_0, a_1, \ldots, a_n)$, denoted by
$\alpha(I(a_0, a_1, \ldots, a_n))$ or $\alpha(I)$ for short, is defined to be the
largest monomial with a nonzero coefficient with respect to the above order. For example, for the semi-invariants
$I_1$ and $I_2$ in \eqref{I1I2} and \eqref{I1I2-b}, we have
\[ \alpha(I_1)= a_1^2 a_2^2, \quad \mbox{and} \quad \alpha(I_2)= a_0 a_2^3. \]

Since the set of semi-invariants of a binary $n$-form of degree $k$ and weight $m$ forms
a vector space, by a triangulation process or the Gauss elimination we may always find a set of semi-invariants whose leading terms are strictly decreasing with respect to the anti-lexicographic order.

We are now in a position to complete the proof of Theorem \ref{strict}.

\noindent
{\it Proof of Theorem \ref{strict}.} Let \[ t=\delta(k-r, n, m-nr/2).\]
If $t=0$,  by the strict unimodality of the Gaussian coefficients established by Pak and Panova \cite{PP13a,PP13b} as stated in Theorem \ref{1}, we find that $\delta(k,n,m)\geq 1$.
If $t>0$,  we assume that $I_1, I_2, \ldots, I_t$, with the variables $a_0, a_1, \ldots, a_n$ being
suppressed, are linearly independent semi-invariants of degree $k-r$ and weight $m-nr/2$. By the triangulation process, we may further assume that the leading terms of
$I_1, I_2, \ldots, I_t$ are strictly decreasing with
 respect to the anti-lexicographic order, that is,
\begin{equation} \label{alpha1t}
 \alpha(I_1) > \alpha(I_2) > \cdots > \alpha(I_t).
 \end{equation}

Next we attempt to construct
$t+1$ linearly independent semi-invariants of degree $k$ and weight $m$. By Lemma \ref{nr8},
 there exist two linearly independent semi-invariants $J_1$ and $J_2$
 of a binary $n$-form of degree $r$  and weight $nr/2$. Without loss of
 generality, let us assume that
 \begin{equation} \label{J12}
 \alpha(J_1) > \alpha(J_2).
 \end{equation}

We claim that $J_1 I_1, J_1I_2, \ldots , J_1I_t, J_2I_t$ are linearly
independent semi-invariants of degree $k$ and weight $m$.
The degree and weight conditions are easily satisfied.
It remains to verify that $J_1 I_1, J_1I_2, \ldots , J_1I_t, J_2I_t$ are linearly independent. To this end, it suffices to show that
the leading terms of $J_1 I_1, J_1I_2, \ldots , J_1I_t, J_2I_t$ are distinct.
For two semi-invariants $K_1$ and $K_2$, it
is evident that
\begin{equation}
\alpha(K_1K_2) = \alpha(K_1) \alpha(K_2).
\end{equation}
Under the assumption \eqref{alpha1t},
we see that
 \begin{equation} \label{alphaIJ}
 \alpha(J_1I_1) > \alpha(J_1I_2) > \cdots > \alpha(J_1I_t).
 \end{equation}
 Invoking the assumption \eqref{J12}, we get
 \begin{equation} \label{alphaJ1J2}
 \alpha(J_1I_t) > \alpha(J_2I_t).
 \end{equation}
Combining \eqref{alphaIJ} and \eqref{alphaJ1J2}, we reach the
conclusion that $J_1 I_1, J_1I_2, \ldots , J_1I_t, J_2I_t$ have distinct leading terms,
and hence they must be linearly independent. This completes the proof.
\qed

 In fact, there is a slightly more general construction of semi-invariants, as stated below. The proof is omitted.

\begin{lem}\label{-1}
Let $k_1,k_2,n\geq 0$, $0\leq m_1\leq nk_1/2$, $0\leq m_2\leq nk_2/2$. For a binary $n$-form,
assume that there are $t_1$  linearly independent semi-invariants of degree $k_1$ and weight $m_1$, and $t_2$ linearly independent semi-invariants of degree $k_2$ and weight $m_2$, where $t_1, t_2\geq 1$. Then there exist $t_1+t_2-1$ linearly independent semi-invariants of degree $k_1+k_2$ and weight $m_1+m_2$.
\end{lem}

We conclude with a remark of Professor S. S. Chern.
He once commented  that Gian-Carlo Rota was fond of
invariant theory, which he considered
as a subject of old mathematics, in his own words, but not in a
negative tone, to avoid any misunderstanding. In the same context,
Professor Chern also mentioned that some papers of Gauss
remain to be explored, but perhaps not easy to comprehend.
As time passes, it may be witnessed that Rota had his reasons.
At least it is our belief that this is something that should not
completely slip our minds, or ``there is something to it'',
as Rota would have put it.

\vspace{0.5cm}
\noindent{\bf\large Acknowledgment.}\ \
This work was done under the auspices of the National Science Foundation of China.

\end{document}